\documentclass[a4paper,10pt,reqno]{amsart}
\usepackage[TS1,T1]{fontenc}

\usepackage[varg]{newtxmath}
\usepackage{tgtermes}
\usepackage[scale=0.92]{tgheros} 
\usepackage{enumerate}

\usepackage[hyphens]{url}
\usepackage[colorlinks,citecolor=blue,urlcolor=blue,linkcolor=blue,breaklinks,bookmarks]{hyperref}
\usepackage[hyphenbreaks]{breakurl}  %

\usepackage[elide]{natbib}

\theoremstyle{remark}

\newcommand{\PP}{\sqsubset}
\newcommand{\nPP}{\nsqsubset}
\newcommand{\Ing}{\sqsubseteq}
\newcommand{\nIng}{\nsqsubseteq}
\newcommand{\Ks}{\exists}
\newcommand{\Ko}{\forall}
\newcommand{\nz}{\mathrm}
\let\AND=\wedge
\let\OR=\vee
\let\Row=\Leftrightarrow
\DeclareMathSymbol{\pExt}{\mathord}{symbolsC}{78}
\DeclareMathSymbol{\Ext}{\mathrel}{symbolsC}{78}

\DeclareMathSymbol{\pOv}{\mathord}{symbolsC}{7}
\DeclareMathSymbol{\Ov}{\mathrel}{symbolsC}{7}
\DeclareMathSymbol{\POv}{\mathrel}{AMSa}{71}

\newcommand{\Sum}{\mathrel{\mathsf{sum}}}
\newcommand{\Sup}{\mathrel{\mathsf{sup}}}

\newcommand{\IMP}{\mathrel{\Longrightarrow}}
\newcommand{\Imp}{\mathrel{\Rightarrow}}

\let\bprod=\sqcap
\let\bsum=\sqcup

\MakeRobust{\ref}

\makeatletter
\newcommand{\labeltext}[2]{%
  \@bsphack
  \csname phantomsection\endcsname 
  \def\@currentlabel{#1}{\label{#2}}%
  \@esphack}
\makeatother

\begin{document}

\title{Different Theories of Parts}

\author{Andrzej Pietruszczak}

\address{Department of Logic, Institute of Philosophy, Nicolaus Copernicus University in Toruń}
\curraddr{Stanisława Moniuszki 16/20, 87-100 Toruń, Poland}
\email{pietrusz@umk.pl}
\thanks{It is the English version of the article ``Różne teorie części', \emph{Przegląd Filozoficzny – Nowa Seria}, 32 (1): 125–148 (2023; published 2024; \url{https://doi.org/10.18778/0138-0680.2024.12}). This research was funded in whole or in part by the National Science Centre (NCN), Poland, grant no.\ 2021/43/B/HS1/03187. For the purpose of Open Access, the author has applied a CC-BY public copyright licence to any Author Accepted Manuscript (AAM) version arising from this submission.}



\begin{abstract}
Stanisław Leśniewski's mereology was formulated in a specific way, deviating from standard formalizations. Nowadays, Leśniewski's theory is presented in the form of an elementary theory or translated into the language of the theory of relational structures. In this article, firstly, we look at <<existentially neutral>> theories, in which we do not postulate the existence of any other collective sets than those obtained from the definition and the basic properties of the relation of being a part. We also examine some <<existentially involved>> theories of parts. Among them, there is Grzegorczykian mereology and Leśniewskian mereology. One of the main principles of mereology is the transitivity of the concept of \emph{being a part of}. This property is often questioned in the literature on the subject. In the final part, we present problems related to the transitivity of this concept.
\end{abstract}

\maketitle

\section{Introduction\label{sec1}}

Mereology emerged as a theory of collective sets (or, in modern times, mereological sums). It was constructed by the Polish logician Stanisław Leśniewski (\citeyear{L16}, \citeyear{L91a}, \citeyear{L27}--\citeyear{L31}, \citeyear{L91b}). Collective sets are certain wholes composed of parts and the very notion of \emph{being a collective set} can be defined by means of the relational notion of \emph{being a part of}.\footnote{Leśniewski is not the creator of the notion of a collective set or a collective class. They are discussed by, for example, Whitehead and Russell in their commentaries in \emph{Principia Mathematica} \citeyearpar{WR}. These sets were used by, among others, Whitehead in his considerations of the philosophy of spacetime \citeyearpar[see, {e.g.}][]{W}. Leśniewski gave two formal definitions of collective sets (classes). In his theory, both definitions are equivalent.} Therefore, mereology can be considered a theory of the “relation of parts to the whole” (from the Greek  $\mathord{\mu}\acute{\mathord{\varepsilon}}\!\mathord{\rho}\mathord{o}\mathord{\varsigma}$ means `part').

Leśniewski's mereology was formulated in a specific way, differing from standard formulations. This theory was built on another system of Leśniewski's, which he called ``Ontology''. Nowadays, Leśniewski's theory is presented in the form of a certain elementary theory or by translating it into the language of the theory of relational structures. It can also be analyzed using the so-called plural logic.

Since the etymology of the word `mereology' corresponds to the meaning of the phrase `theory of parts', the former can be used to refer to any formal or semi-formal considerations of parts and not only to Leśniewski's theory. We believe, however, that this may cause some terminological confusion. When considering theories weaker than Leśniewski's theory, we should rather add appropriate qualifying adjectives, as Peter \citet{S87} did, for example, when he examined ``Minimal Extensional Mereology''. In \citep{ja13,ja20}, various existentially neutral and existentially involved theories of parts are analyzed. In the former, we do not postulate the existence of any other collective sets than those that we obtain from the basic properties of the relation of \emph{being a part}.\footnote{Note that collective sets composed of only physical (material) objects are to be objects of the same kind.} For example, in such theories, we will not obtain the existence of a collective set formed from two objects that are parts of a third one. We do not, therefore, postulate the existence of a collective set consisting of the right and left hands of a given person.

Leśniewski's mereology is classified as an existentially involved theory. Their existential involvement consists in the fact that they have additional axioms postulating the existence of collective sets of various groups of objects. Some of such sets can be considered as objects obtained \emph{ad hoc}, which is therefore controversial. For example, it is difficult to accept that there is a material object that would be a collective set consisting of the Moon and the heart of a given person \citep{ja00a,ja00b,ja05,ja18}. What is more, the existence of a separate object that would be a collective set of the right and left hands of a given person is even problematic. In Leśniewski's theory, on the other hand, the unlimited existence of collective sets for all (non-empty) groups of objects (including infinite ones) is postulated. It seems that such a solution is only admissible in pointless geometry and pointless topology, which concern spatial regions or space-time events.\footnote{Let us add that in these pointless theories, we have points. However, they are not assumed as primitive. They are defined on the basis of primitive concepts such as spheres, solids, or regions, using the relation of \emph{being of a part of} \citep[see, {e.g.},][]{T29,T56,G60, GP08,GP09,GP18a, GP18b, GP19,GP23}.}

In Section~\ref{sec2}, we will present the basic concepts of mereology. In Section~\ref{sec3}, we will deal with mereological sums as collective classes of a given group of objects. In Section~\ref{sec4}, we will present existentially neutral theories. Theories that are existentially involved will be presented in Section~\ref{sec5}. Among them will be the strongest of them---Leśniewskian mereology and two theories proposed by Andrzej \citet{G55}. In the final section, we will sketch a problem related to the transitivity of the concept of \emph{being a part of}. This property is often called into question in the literature. We present an analysis without the assumed transitivity.

\section{Basic concepts of mereology\label{sec2}}

In this section, we will present the basic concepts of the theory of parthood and their fundamental formal properties. We will assume that the relational concept of \emph{being a part of} should create a strict partial order in any universe of considerations, i.e., it should be transitive and irreflexive, which also gives asymmetry. Moreover, in all non-degenerate universes (i.e. those that have at least two elements), there should be no smallest object (zero), and instead, there should be two objects that have no common part.

\subsection{Parts as fragments\label{sec:czjk}}
In everyday speech, the expression `part' is usually understood as having the sense of the expressions `fragment', `bit' and `piece' when we relate them to spatial objects (regions) or spatiotemporal events. Thus understood, the relation of a part to the whole has two properties:\vspace{-3pt}
\begin{enumerate}[(a)]
\item\labeltext{a}{c-a} no object is its own part;
\item\labeltext{b}{c-b} there are not two objects such that the first could be a part of the second and the second a part of the first.
\end{enumerate}
Thanks to condition \eqref{c-a}, we have no difficulty in interpreting the phrase `two objects' in condition~\eqref{c-b}. It can be seen that it concerns ``two different'' objects. Sentence \eqref{c-a} says that the relation of part to whole is irreflexive and \eqref{c-b} says the this concept is antisymmetric. Both conditions taken together are equivalent to the fact that this concept is asymmetric. In order to abbreviate these and other properties of the concept \emph{being a part of}, we will symbolise `$x$ is a part of $y$' with `$x \PP y$' and likewise with other variables. We are therefore treating the symbol `$\PP$' as a two-argument predicate. With respect to an arbitrary universe of discourse $U$ (this being a non-empty distributive set), the irreflexivity antisymmetry and asymmetry of the concept \emph{being a part of} may be expressed in the form of the following two \pagebreak conditions:
\begin{gather*}
\neg\Ks_{x\in U}\; x \PP x, \label{pz-PP} \tag{$\nz{irr}_{\scriptscriptstyle\PP}$}\\
\neg\Ks_{x,y\in U}\bigl(x\neq y\AND x\PP y \AND y\PP x\bigr), \label{antys-PP} \tag{$\nz{antis}_{\scriptscriptstyle\PP}$}\\
\neg\Ks_{x,y\in U}\bigl(x\PP y \AND y\PP x\bigr). \label{as-PP} \tag{$\nz{as}_{\scriptscriptstyle\PP}$}
\end{gather*}
The conjunction of \eqref{pz-PP} and \eqref{antys-PP} is logically equivalent to \eqref{as-PP}.

Leśniewski assumed that the relation of part to whole is asymmetric (i.e.\ also irreflexive and antisymmetric) and transitive, i.e.
\begin{enumerate}[(c)]
\item\labeltext{c}{c-c} every part of some part of a given object is also a part of the given object.
\end{enumerate}
Therefore, the following condition must be met:
\begin{equation}\label{p-PP}\tag{$\nz{t}_{\scriptscriptstyle\PP}$}
\Ko_{x,y,z\in U}\bigl((x\PP y \AND y\PP z) \IMP x\PP z\bigr).
\end{equation}
\indent
In support of the transitivity of the concept of \emph{being a part}, the following example is often given: my left arm is a part of my body, from which it follows that my left hand is part of my body. Nicholas \citet{R55}, however, shows that the transitivity of the relation of part to whole is in other cases problematic. He provides the following counterexample: a nucleus is part of a cell and a cell is part of an organ, but the nucleus is not part of an organ. This is so, at least, if we consider a part to be a direct functional constituent of a whole, for a nucleus is not a part of an organ.  \citet[pp.~107--108]{S87} has observed, however, that the concept of a part understood transitively corresponds to spatiotemporal inclusion and in that sense a nucleus is part of an organ. Simons claims that the fact that the word `part' has an additional meaning does not undermine the mereological concept of a part, because there is no reason to think that the mereological concept should embrace all the meanings of the word `part' but only those that are basic and most important. The transitivity of the relation \emph{is a part of} does not therefore present special difficulties when referring to spatiotemporal relations, including events.

From the irreflexivity and transitivity of the part-whole relation, it follows that it is acyclic, i.e., there are no closed cycles regarding being a part. The following scheme for any positive natural number $n$ expresses this:
\[\label{ac-PP}\tag{$\nz{ac}_{\scriptscriptstyle\PP}$}
\neg\Ks_{x_1,\ldots,x_n\in U}\bigl(x_1\PP x_n \AND\dots\AND x_{n}\PP x_1\bigr).
\]
Of course, the above scheme also expresses \eqref{pz-PP} and \eqref{as-PP} (for $n= 1, 2$).

If we reject the transitivity of the concept of\emph{ being a part of}, we assume its acyclicity, which gives asymmetry and irreflexivity. In the final section, we will sketch the problems related to the transitivity of this concept.

\subsection{Another meaning of the word `part'}

In the literature on the subject, the custom has spread of using the expression `proper part' rather than `part'. The term `part' here takes on a new meaning, in which it has wider range of use. A~part of a given object is now understood to be the object itself or any given part in the ordinary sense of that term. A~part of a given object that differs from that object itself is named a \emph{proper part}. Therefore, in conditions \eqref{c-a} and \eqref{c-b} the term `proper part' would appear in place of `part'. It follows directly from the new meaning of `part' that \emph{every object is its own (improper) part}. If we understand `two objects' in the sense of ``two different'' objects, then, under the new meaning of `part', condition \eqref{c-b} is satisfied: that is, we have the antisymmetry of the altered concept \emph{being a part}.

If we use a given word in a new meaning, we must treat it as a new concept and apply a new symbolic designation to it. Thus, with this new meaning, the phrase `$x$~is part of~$y$' will be written as `$x\Ing y$'. The relationship between the two concepts is expressed by the following formula, defining for all $x,y\in U$, the new concept in terms of the old one:
\[\label{def-Img}\tag{$\nz{df}\,{\Ing}$}
x\Ing y \iff (x=y \OR x\PP y).
\]
The reflexivity of the predicate `$\Ing$'  follows directly from the reflexivity of the identity predicate `$=$'. Moreover, the antisymmetry of `$\Ing$'is obtained from \eqref{antys-PP} and the identity property. So we obtain:\vspace{-6pt}
\begin{gather*}
\Ko_{x\in U}\; x \Ing x\, , \label{z-Ing}\tag{$\nz{r}_{\scriptscriptstyle\Ing}$}\\
\neg\Ks_{x,y\in U}\bigl(x\neq y\AND x\Ing y \AND y\Ing x\bigr). \label{antys-Ing} \tag{$\nz{antis}_{\scriptscriptstyle\Ing}$}
\end{gather*}
Similarly, if we assume that the relation $\PP$ is transitive, then so is the relation $\Ing$, i.e.:
\begin{equation}\label{p-Ing}\tag{$\nz{t}_{\scriptscriptstyle\Ing}$} \Ko_{x,y,z\in U}\bigl((x\Ing y \AND y\Ing z) \IMP x\Ing  z\bigr).
\end{equation}
Moreover, by \eqref{pz-PP} and \eqref{as-PP}, for all $x$ and $y$ we obtain:
\begin{equation*}
\begin{split}
x\PP y &\iff (x\Ing y\AND x\neq y),\\
&\iff (x\Ing y\AND y\nIng x).
\end{split}
\end{equation*}
The above formulas are not definitions of the relation $\PP$ because it is primitive in our case. The broadening of the extension of the word `part' can sometimes lead to ``philosophical misunderstandings''.

\subsection{Ingredienses} Leśniewski did not change the ordinary meaning of the word `part'. Instead, he adopted the word `ingredjens' in his works during the years \citeyear{L27}--\citeyear{L31}, a~word which did not exist in Polish at that time. It seems that he felt the need for a new term. We will also employ his neologism but spell it thus: `ingrediens'. The strange-sounding term `ingrediens' is in this case an <<ally>>, as it reminds us that it is an <<artificial concept>>.\footnote{Let us note that a misunderstanding arises from the choice of `ingredient' as the translation of `ingredjens' in the English translations of Leśniewski's works \citep[{cf.}][]{L91b}. The word `ingredient' has a ready translation in Polish as `\emph{składnik}', and Leśniewski deliberately distanced himself from this word (since the whole is not its ingredient). It is, therefore, better to use `ingrediens' (`ingredienses’) in translation, thus preserving continuity with its use in Polish Leśniewski's works.} An ingrediens of a given object is therefore it itself and all its parts (in the ordinary sense of the word `part').

\subsection{The absence of an <<empty (zero) element>>}
When we assume that the universe of discourse consists of physical objects or spatial regions or space-time events, respectively, we assume that there is more than one such object. Furthermore, we exclude from our considerations the existence of an <<empty object>>, an <<empty region>>, or an <<empty event>>, which would be a part of any other object, region or event, respectively. We, therefore, have no analogy to set theory---the theory of distributive sets (classes)---where we assume the existence of an empty set $\emptyset$, which is a subset of every distributive set (class).

In an algebraic sense, such an <<empty element>> would correspond to the \emph{zero}, i.e., the least element with respect to~$\Ing$. We say that
\begin{enumerate}[\textbullet]
\item $x$ is the \emph{zero} if and only if $\Ko_{u\in U}\; x \Ing u$.
\end{enumerate}
We will allow the existence of the zero when and only when the universe has exactly one element, which is simultaneously the zero and the greatest element with respect to~$\Ing$, which we call the \emph{unity}. We say that
\begin{enumerate}[\textbullet]
\item $x$ is the \emph{unity} if and only if $\Ko_{u\in U}\; u\Ing x$.
\end{enumerate}
We do not have to note that such element which simultaneously is the zero and the unity is <<empty>>. It is not an interesting case because we are dealing with a degenerate structure.

In all the theories considered, there is no zero in any non-degenerate universe:
\begin{equation*}\label{nepe}\tag{$\nexists 0$}
\Ks_{z,u\in U}\; z\neq u \IMP \neg\Ks_{x\in U}\; x \text{~is zero}.
\end{equation*}
Principle \eqref{nepe} will follow from other principles to be later adopted. We give it here because it allows the further auxiliary concepts to take on their proper meaning.

\subsection{The relation of \emph{being exterior to}} The second auxiliary binary relation which may hold between elements of the universe $U$ is the relational concept of \emph{is exterior to}, which we will denote by `$\pExt$'. We say that one object is exterior to another if they have no common ingredients. In symbolic notation, for all $x,y\in U$, we put:
\begin{equation}\label{def-Ext}\tag{$\nz{df}\,\pExt$}
x\Ext y \iff \neg\Ks_{z\in U}(z\Ing x \AND z\Ing y).
\end{equation}
The name and definition of this relation are due to Leśniewski. Of course, $\Ext$ is symmetric, so we may say that ``objects $x$ and $y$ are exterior with respect to one another''. Since $\Ing$ is reflexive, $\Ext$ is irreflexive. Moreover, the intuitions related to the meaning of the phrase `is exterior to' are taken from the case when $\PP$ is a <<true>> relation of \emph{being a part of}, and $\Ing$ is a <<true>> relation of \emph{being an ingredient of}, and remembering that no element is zero in a non-degenerate universe (see \eqref{nepe}).\footnote{The fact that $x\Ext y$ holds does not exclude the fact that $x$ and $y$ are in contact. The only thing is that $x$ and $y$ do not have any part in common. The distinction between the tangent and non-tangent cases of being exterior is possible only in point-free geometry or topology \citep[see, e.g.,][]{T29, T56, G60, GP08, GP09, GP18a, GP18b, GP19, GP23}. Thus, it is completely different than in ordinary point geometry, where tangent figures are not disjoint because they have a common point.}

That two objects are exterior to each other is equivalent to the fact that neither of them is part of the other and that they have no part in common, i.e., for all $x,y\in U$:
\begin{equation*}
\begin{split}
x\Ext y &\iff \bigl(x\ne y\AND x\nPP y \AND y\nPP y\AND \neg\Ks_{z\in U}(z\PP x\AND z\PP y) \bigr),\\
&\iff \bigl(x\nIng y \AND y\nIng y\AND \neg\Ks_{z\in U}(z\PP x\AND z\PP y) \bigr)
\end{split}
\end{equation*}
\indent Thanks to principle \eqref{nepe}, the relation $\pExt$ does not become <<automatically>> empty in non-degenerate structures. This said, it is not guaranteed that it won't be empty. To determine  whether it is, we must consider those structures that satisfy the following principle (which is stronger than \eqref{nepe}):
\begin{equation}\label{EExt}\tag{$\Ks\pExt$}
\Ks_{z,u\in U}\; z\neq u \IMP \Ks_{x,y\in U}\; x\Ext y   .
\end{equation}
\subsection{The relations of \emph{overlapping} and \emph{crossing}}
The third auxiliary binary relation which can hold between elements of the universe $U$ is the relational concept of \emph{overlapping}, which we will denote by `$\pOv$'. Its designation and name come from the word `overlap', used in \citep[p.~47]{LG}. We say that two objects overlap when they have at least one common ingredient. In symbolic notation, for all $x,y\in U$, we put:
\begin{equation}\label{def-Ov}\tag{$\nz{df}\,\pOv$}
x\Ov y \iff \Ks_{z\in U}(z\Ing x \AND z\Ing y).
\end{equation}
Of course, $\pOv$ is symmetric, so we may say that ``objects $x$ and $y$ overlap''. Since $\Ing$ is reflexive, $\Ext$ is also reflexive. Furthermore, the relations $\pOv$ and $\pExt$ complement each other. Thanks to principles \eqref{nepe} and \eqref{EExt}, our intuitive understanding of the relation should become clearer. Note, however, that the condition $x\Ov y$ does not mean that $x$ and $y$ cross since we do not exclude that either are identical or one of them is part of the other (which conflicts with the meaning of the word `overlap'). Namely, the fact that $x\Ov y$ is equivalent to the fact that either $x=y$, or $x$ is part of $y$, or vice versa, or $x$ and $y$ have a common part. Formally:
\begin{equation*}
\begin{split}
x\Ov y&\iff x= y\OR x\PP y\OR y\PP x\OR\exists_{z\in U}(z\PP x\AND z\PP y)\bigr), \\
&\iff x\Ing y \vee y\Ing x\vee \Ks_{z\in U}(z\PP x\AND z\PP y).
\end{split}
\end{equation*}
\indent The word `overlap' is more closely related to the relational concept of \emph{overlapping}, which we will denote by `$\pOv$'.relation of \emph{crossing}, which we will denote with the symbol `$\POv$'. We say that two (different) objects cross one other (or properly overlap one another) iff they have some part in common, but neither of them is part of the other (i.e., that common part is not identical to either of them). Formally, we introduce into $U$ the binary relation $\POv$ via the following definition:
\begin{equation}\label{def-POv}\tag{$\nz{df}\,\mathord{\POv}$}
x \POv y \iff x\neq y\AND x\nPP y \AND y\nPP x \AND \Ks_{z\in U}(z\PP x\AND z\PP y).
\end{equation}
This relation is symmetric and irreflexive. We will also express it as follows:
\begin{equation*}
\Ko_{x,y\in U}\bigl(x \POv y \iff (x\Ov y \AND x\nIng y\AND y\nIng x)\bigr).
\end{equation*}
Therefore, two objects cross one other if and only if they overlap one another, but neither is an ingrediens of the other. Thanks to the principles \eqref{nepe} and \eqref{EExt}, the intuitions related to the relations $\pOv$ and $\POv$  become clear. Moreover, by virtue of the adopted definitions, \eqref{pz-PP} and \eqref{as-PP}, we have the following connections between the introduced relations:
\begin{align*}
\Ko_{x,y\in U}\bigl( x\Ov y &\iff x= y\OR x\PP y\OR y\PP x\OR x\POv y\bigr),\\
\Ko_{x,y\in U}\bigl(x\nIng y&\iff x\Ext y\OR x\POv y \OR y\PP x\bigr) .
\end{align*}
Thus, two (different) objects overlap if and only if either they cross or one is part of the other. Furthermore, one object is not an ingrediens of the other if and only if either they are exterior to each other, or they cross, or the other is part of the first.
\section{Mereological sums and supprema\label{sec3}}

\subsection{Collective classes as mereological sums}
Let `P' be a schematic letter representing a general name of certain elements of the universe of discourse $U$. The extension of this name is the distributive set of Ps belonging to $U$, i.e., the set $\{u\in U: u$ is a P\}. Following Alfred \citet{T29, T56}, instead of writing that $x$ is a collective class of all Ps, we will write that $x$ is a mereological sum of all elements of the (distributive) set of all Ps (short: the mereological sum of all Ps).  We can also go straight to the consideration of arbitrary subsets of the universe without restricting ourselves to the extensions of general names. It allows us to leave schematic letters out of our analysis. Instead of a mereological sum of all Ps, we can speak of the sum of all elements of a given subset of the set $U$. For an arbitrary element $x$ from $U$ and arbitrary subset $S$ of $U$, Leśniewski's definitional schema, featuring the expression ``mereological sum'', may be written as follows:

\begin{center}
\parbox{114mm}{$x$ is a \emph{mereological sum} of all elements of a set $S$ if and only if the following two conditions are satisfied:\vspace{-3pt}
\begin{enumerate}[\textbullet]
\item every element of the set $S$ is an ingrediens of $x$,
\item every ingrediens of $x$ overlaps some element of~$S$.
\end{enumerate}}
\end{center}

\noindent Obviously, we could have written the above definition with the help of the primitive concept of \emph{being a (proper) part} instead of the concepts of \emph{being an ingrediens} and \emph{overlapping}.

However, instead of\vspace{-3pt}
\[
\text{$x$ is a mereological sum of all elements of a set $S$}
\]
we can write:\vspace{-6pt}
\[
x\Sum S   .
\]
To put things formally, we are introducing $\Sum$ as a binary relation included in the Cartesian product $U\times \mathcal{P}(U)$. We define this relation according to the recipe above by putting for arbitrary $x$ from $U$ and $S$ from $\mathcal{P}(U)$:
\begin{equation}\label{def-Sum}\tag{$\nz{df}\,\mathord{\Sum}$}
x\Sum S \iff\Ko_{s\in S}\; s\Ing x \AND \Ko_{u\in U}(u\Ing x\Imp \Ks_{s\in S}\; s\Ov u).
\end{equation}
As a consequence of the reflexivity of $\Ing$, for arbitrary $x\in U$ and $S\in\mathcal{P}(U)$, we get:
\[
x \Sum S \IMP S\neq \emptyset,
\]
i.e., there is no mereological sum of the empty set $\emptyset$. It says that there is no such thing as an empty collective class. It agrees with \eqref{nepe}, stating that in a non-degenerate universe, no zero element plays the role of an empty object. Moreover, it is consistent with the following general rule, which we get from the reflexivity of the relation $\pOv$:\vspace{-3pt}
\begin{enumerate}[\textbullet]
\item if $x\in S$ and $S$ is included in the set of all ingredienses of $x$, then $x\Sum S$, i.e.,\\ if $x\in S \subseteq \{u\in U : u\Ing x\}$, then $x\Sum S$.
\end{enumerate}
Hence we get:
\begin{enumerate}[\textbullet]
\item $x$ is a mereological sum of itself, i.e.\ $x\Sum\{x\}$;
\item $x$ is a mereological sum of all its ingredienses, i.e.\ $x\Sum\{u\in U:u\Ing x\}$.
\end{enumerate}

Further, we get that there are no objects with exactly one part (see p.~\pageref{eq:nie-ma-jednej}), although there may be objects that have no parts at all---so-called \emph{mereological atoms}. We obtain:
\begin{enumerate}[\textbullet]
\item if $x$ is not an atom, then $x$ is a mereological sum of its parts, i.e.\ $x\Sum\{u\in U:u\PP x\}$.
\end{enumerate}

From \eqref{antys-Ing}, we get that a given universe $U$ has at most one unity.  So far, we have not dealt with the problem of whether in $U$ there is the unity. This problem is related to the existence of the mereological sum of the entire universe. Namely, we have:
\begin{enumerate}[\textbullet]
\item $x\Sum U$ if and only if $x$ is the unity in $U$.
\end{enumerate}
The universe can, therefore, have at most one mereological sum (if it has unity).
\subsection{The uniqueness of mereological sum}
According to its meaning, the phrase `collective class of all Ps' can have at most one referent. Such an axiom was accepted by \citet{L27,L28,L91b}. In our terminology, the counterpart of his axiom says that if a (distributive) set has a mereological sum, it is unique:
\begin{equation}\label{fun-Sum}\tag{$\nz{U}_{\scriptscriptstyle\Sum}$}
\Ko_{S\in\mathcal{P}(U)}\Ko_{x,y\in U}\bigl((x\Sum S\AND y\Sum S)\IMP x=y\bigr).
\end{equation}
Previously assumptions only allowed us to say that if the universe has a sum, then only one (its unity if it exists). Now, every distributive set can have at most one sum. So, $\{x\}$ is the only singleton whose sum is $x$ (since  for any $y$, we have $y\Sum\{y\}$):
\[\label{sf-Sum}\tag{$\nz{S}_{\scriptscriptstyle\Sum}$}
\Ko_{x,y\in U}(x\Sum \{y\}\IMP x=y).
\]
Therefore, (since $x\Sum\{x\}$) we can identify the mereological sum of a singleton $\{x\}$ with $x$ itself, which we write as $x=\llbracket x \rrbracket$. Therefore, the collective class constructed from one object is that object. We do not claim that it has exactly one element because each part of the object $x$ must be considered a mereological element of the collective class $\llbracket x\rrbracket$.

Finally, note that \eqref{fun-Sum} entails the following principle of extensionality with respect to~$\PP$:
\[\label{ext-PP}\tag{$\nz{ext}_{\PP}$}
\Ko_{x,y\in U}\big((\Ks_{z\in U}\; z\PP x\AND \Ko_{u\in U}(u\PP x\Row u\PP y )) \IMP x=y\big).
\]
Indeed, let $P_x$ and $P_y$ be the sets of all parts of $x$ and $y$. By the assumption, $P_x=P_y\neq\emptyset$. We know that $x\Sum P_x$ and $y\Sum P_y$. So $x=y$, by \eqref{fun-Sum}. It is necessary to assume that $x$ is not an atom. The converse implication cannot be obtained since $x$ may have no part. There are, therefore, no two objects which have parts at all and which are the same parts.

The principle of extensionality with respect to $\Ing$ follows from \eqref{z-Ing} and \eqref{antys-Ing}:
\[\label{ext-Ing}\tag{$\nz{ext}_{\Ing}$}
\Ko_{x,y\in U}\big(\Ko_{u\in U}(u\Ing x\Row u\Ing y) \IMP x=y\big).
\]
Indeed, from \eqref{z-Ing} we get: $\Ko_{u\in U}(u\Ing x\Imp u\Ing y)\Imp x \Ing y$ and $\Ko_{u\in U}(u\Ing y\Imp u\Ing x)\Imp y\Ing x$. Hence and from \eqref{antys-Ing} we have \eqref{ext-Ing}.
\subsection{Suprema} The universes under consideration are partially ordered by the relation $\Ing$, which is reflexive, antisymmetric, and transitive in them. Thus, we can consider the binary relation of \emph{being a least upper bound} (or \emph{supremum}) in them. A least upper bound of a subset $S$ of a universe $U$ is to be the smallest object in $U$ such that all members of $S$ are its ingrediens. The fact that an object $x$ is a least upper bound (supremum) of a set $S$ will be written as $x \Sup S$. We can, therefore, construct a binary relation $\Sup$ contained in the set $U\times \mathcal{P}(U)$, assuming that for all $S\in\mathcal{P}(U)$ and $x\in U$:
\begin{equation}\label{def-Sup}\tag{$\nz{df}\,\mathord{\Sup}$}
x\Sup S \iff\Ko_{s\in S}\; s\Ing x \AND \Ko_{u\in U}(\Ko_{s\in S}\; s\Ing u \Imp x\Ing u).
\end{equation}
\indent By applying \eqref{z-Ing} we obtain: $x\Sup\{x\}$ and $x\Sup\{u\in U : u\Ing x\}$. Moreover, from \eqref{antys-Ing} we get the uniqueness of the relation $\Sup$, i.e.:
\begin{equation}\label{fun-Sup}\tag{$\nz{U}_{\scriptscriptstyle\Sup}$}
\Ko_{S\in\mathcal{P}(U)}\Ko_{x,y\in U}\bigl(x\Sup S \AND y\Sup S\IMP x=y \bigr) ,
\end{equation}
and from \eqref{z-Ing} and \eqref{antys-Ing} we have: $\Ko_{x,y\in U}(x\Sup\{y\}\IMP x=y)$.

From the assumptions made so far, it is not possible to derive any interesting connections between the relations $\Sum$ and $\Sup$. We have introduced this second relation only in order to show the connections between these two relations when making further assumptions.

\section{Existentially neutral theories of parts\label{sec4}}

As we mentioned in the introduction, in existentially neutral theories, the existence of mereological sums follows only from the accepted definitions and basic properties of the part-relation. We do not have any additional axioms postulating the existence of such sums. Let us add that such an existential axiom does not have to postulate the existence of a mereological sum explicitly.

\subsection{The Weak Supplementation Principle}
Examples of neutral existential axioms are Simons' \citeyearpar{S87} two supplementation principles: weak and strong. The first one (\emph{Weak Supplementation Principle}) has the form:
\begin{equation}\label{WSP}\tag{$\nz{WSP}$}
\Ko_{x,y\in U}\bigl(y\PP x \IMP \Ks_{z\in U}(z\PP x\AND z\Ext y)\bigr),
\end{equation}
i.e. if one object is a part of another, then some object is a part of the second and is exterior to the first. From \eqref{WSP} and the irreflexivity of the relation $\Ext$, it follows that no object has exactly one part, i.e.:
\begin{equation*}\label{eq:nie-ma-jednej}
\Ko_{x,y\in U}\bigl(y\PP x\IMP\Ks_{z\in U}(z\PP x\AND z\ne y)\bigr),
\end{equation*}
\indent In \citep{ja00b,ja18} it is shown that \eqref{WSP} entails \eqref{p-PP}. Indeed, assume that $x\PP x$. Then, by \eqref{WSP}, for some $z$, we have: $z\PP x$ and $z\Ext x$. This first gives $z\Ov x$, so we obtain a contradiction. Therefore, given \eqref{WSP} and \eqref{p-PP}, one need not assume either \eqref{pz-PP} or \eqref{as-PP}.\footnote{\citet{S87}, and following him other authors assume \eqref{WSP}, \eqref{p-PP} plus either \eqref{pz-PP} or \eqref{as-PP}.} \eqref{WSP} entails also both principles \eqref{nepe} and \eqref{EExt}.

It has been proven that \citep[see][pp.~71--72, 83--84, 59–60, 47, respectively]{ja00b,ja18,ja13,ja20}:
\begin{enumerate}[\textbullet]
\item \eqref{WSP} is equivalent to the conjunction \eqref{pz-PP} and \eqref{sf-Sum}.
\end{enumerate}
Hence we get that:
\begin{enumerate}[\textbullet]
\item \eqref{pz-PP} and \eqref{fun-Sum} entail \eqref{WSP}.
\end{enumerate}

We have already mentioned that in order to obtain any connection between the relations $\Sum$ and $\Sup$, we need to assume something additional about the relation $\PP$. The first such assumption is \eqref{WSP}. In \citep[pp.~147, 180, 63, 52, respectively]{ja00b,ja18,ja13,ja20}, it was shown that \eqref{WSP} is equivalent to the following sentence:
\begin{equation}\label{Sum-Sup-id}\tag{$\diamond$}
\Ko_{S\in\mathcal{P}(U)}\Ko_{x,y\in U}\bigl(x\Sum S \AND y\Sup S\IMP x=y \bigr) ,
\end{equation}
which says that if the sum and the supremum of a given set exist, then they are equal.

\subsection{The Strong Supplementation Principle\label{4.2}}
The second principle adopted by Simons (1987)---the \emph{Strong Supplementation Principle}---states that if one object is not an ingrediens of another, then some object is an ingrediens of the first and is exterior to the second:
\begin{equation}\label{SSP}\tag{$\nz{SSP}$}\
\Ko_{x,y\in U}\bigl(x\nIng y\IMP\Ks_{z\in U}(z\Ing x\AND z\Ext y)\bigr) .
\end{equation}
The successor must be $\Ing$, not $\PP$ because $x$ can be an atom.

Notice that \eqref{as-PP} and \eqref{SSP} entail \eqref{WSP}. Indeed, let $y\PP x$. Then, by \eqref{pz-PP} and \eqref{as-PP}, we have $x\nIng y$. Hence, by \eqref{SSP}, for some $z$, we have: $z\Ing x$ and $z\Ext y$. Hence, by the assumptions, we get $z\neq x$. So $z\PP x$.

So we see that under the standard assumptions for the relation $\PP$, \eqref{SSP} is stronger than \eqref{WSP}. We previously showed that \eqref{fun-Sum} is also stronger than \eqref{WSP}. It turns out that in terms of forcer, \eqref{fun-Sum} is between \eqref{SSP} and \eqref{WSP}. In \citep[pp.~364, 220--221, respectively]{ja00a, ja05}, it was shown that \eqref{antys-Ing}, \eqref{p-Ing} and \eqref{SSP} entail \eqref{fun-Sum}. Indeed, from \eqref{p-Ing} and \eqref{SSP}, for all $x,y\in U$ and $S,Z\in \mathcal{P}(U)$, we have:
\[\label{ast}\tag{$\ast$}
\bigl(\Ko_{u\in U}(u\Ing x\Ing \Ks_{s\in S}\; s\Ov u)\AND S\subseteq Z\wedge \Ko_{z\in Z}\; z\Ing y\bigr)\IMP x\Ing y.
\]
Hence, taking $S=Z$, we get that if $x\Sum S$ and $y\Sum S$, then $x\Ing y$ and $y\Ing x$, i.e., $x=y$, by \eqref{antys-Ing}.

From \eqref{SSP} the following version of the supplementation principle follows:
\[
\Ko_{x,y\in U}\bigl(x\POv y\IMP \Ks_{z\in}(z\PP x\AND  z\Ext y)\bigr),
\]
i.e., if two objects cross, then one of them has an exterior part to the other. Indeed, if $x\POv y$, then $x\Ov y$ and $x\nIng y$. Hence, by \eqref{SSP}, for some $z$, we have $z\Ing x$ and $z\Ext y$. Hence, we have $z\neq x$ since $x\Ov y$. So $z\PP x$.

In \citep[pp.~74--76, 90, 47--48, 35, respectively]{ja00b,ja18,ja13,ja20} it was proved that \eqref{SSP} is equivalent to each of the following two principles:
\begin{align*}
\Ko_{x,y\in U}\bigl(\Ko_{u\in U}(u\Ov x\Imp u\Ov y)&\IMP x\Ing y\bigr), \label{SSPov}\tag{$\nz{SSP}_{\scriptscriptstyle\Ov}$}\\
\Ko_{x,y\in U}\bigl(\Ko_{u\in U}(u\Ext y\Imp u\Ext x)&\IMP x\Ing y\bigr). \label{SSPext}\tag{$\nz{SSP}_{\scriptscriptstyle\Ext}$}
\end{align*}
The equivalence of the two above principles is obtained from the relation between $\Ov$ and $\Ext$. From \eqref{z-Ing} we have the converse implications to \eqref{SSPov} and \eqref{SSPext}, which together gives:
\begin{align*}
\Ko_{x,y\in U}\bigl(x\Ing y&\iff\Ko_{u\in U}(u\Ov x\Imp u\Ov y)\bigr), \\ 
\Ko_{x,y\in U}\bigl(x\Ing y&\iff\Ko_{u\in U}(u\Ext y\Imp u\Ext x)\bigr). 
\end{align*}
Thus, one object is an ingrediens of another if and only if every object that overlaps the first overlaps the second, or equivalently if every object exterior to the second is exterior to the first. The relation $\Ing$ is, therefore, expressible in terms of each of the relations $\Ov$ and $\Ext$.

Principle \eqref{SSP} provides another connection between the relations $\Sum$ and $\Sup$. It is stronger than the one provided by \eqref{WSP}. Namely, in \citep{ja00a, ja05}, it is noticed that \eqref{ast} entails the inclusion ${\Sum}\subseteq {\Sup}$, which says that every mereological sum of a given set is also its supremum. So, since only non-empty sets have mereological sums:
\[\label{2ast}\tag{$\ast\ast$}
\Ko_{S\in\mathcal{P}(U)}\Ko_{x\in U}\bigl(x\Sum S \IMP (S\neq\emptyset \AND x\Sup S)\bigr).
\]
The above and \eqref{fun-Sup} entail \eqref{Sum-Sup-id}. Moreover, in \citep[pp.~78, 97--98, respectively]{ja00b, ja18}, it is proved that:
\begin{enumerate}[\textbullet]
\item principle \eqref{SSP} is equivalent to the inclusion ${\Sum}\subseteq {\Sup}$.
\end{enumerate}

Let us note that in existentially neutral theories, we will not obtain either the converse implication to \eqref{2ast} or the inclusion ${\Sup}\subseteq{\Sum}$. We will show later that these last two conditions require the existence of mereological sums, which we will not obtain from the accepted definitions and basic properties of the relation $\PP$.

In addition to the previously discussed definition of the concept of \emph{being a collective class of}, \citet{L31} gave a second explanation of this concept. In Leśniewski's mereology, both definiens of his definitions are equivalent. So, for all $x\in U$ and $S\in \mathcal{P}(U$) we have:\vspace{-3pt}
\begin{align*}
x\Sum S&\iff \Ko_{u\in U}(u\Ext x\Row \Ko_{s\in S}\; s\Ext u), \label{DolarExt}\tag{\text{\$}$_{\scriptscriptstyle\Ext}$}\\
x\Sum S&\iff \Ko_{u\in U}(u\Ov x \Row \Ks_{s\in S}\; s\Ov u).
\label{DolarOv}\tag{\text{\$}$_{\scriptscriptstyle\Ov}$}
\end{align*}
Thus, $x$ is the mereological sum of all elements of $S$ if and only if the fact that a given object is exterior to $x$  (resp.\ overlaps $x$) is equivalent to the fact that this object is exterior to all elements of $S$ (resp.\ overlaps some element of $S$). In the light of the relationship between $\Ext$ and $\Ov$, the above theses are equivalent. Moreover, they follow from \eqref{SSP} itself and their parts ``$\Leftarrow$'' are equivalent to \eqref{SSP}. The proof of these facts can be found in \citep{ja00b, ja05, ja13,ja18,ja20}.

\subsection{The proper parts principle}
Simons adopted it in \citeyearpar{S87}.\footnote{\citet{S87} uses the terminology that the term ` proper part' corresponds to our term `part'.} The following formula expresses it:
\[\label{PPP}\tag{PPP}
\Ko_{x,y\in U}\bigl((\Ko_{z\in U}\;  z\PP x \AND \Ko{u\in U}(u\PP x\Imp u\PP y) ) \IMP x\Ing y\bigr).
\]
Therefore, if every part of a non-atomic object is a part of another, then the former is an ingrediens of the latter. It is necessary to assume that the former is not an atom. Furthermore, $x\PP y$ cannot occur in the consequent of \eqref{PPP} since the antecedent does not entail that $x$ and $y$ are distinct.

Principle \eqref{PPP} follows from \eqref{SSP}. Indeed, let $x$ have a part, and each of its parts is a part of $y$. Then $x\Ov y$. Assume, for the sake of contradiction, that $x\nIng y$. Then, by \eqref{SSP}, for some $z$, we have $z\Ing x$ and $z\Ext y$. Hence $z\neq x$ since $x\Ov y$. Therefore, $z\PP x$. Hence, by the assumption, also $z\PP y$. So, we have a contradiction.

\subsection{Other extensionality principles}
By \eqref{antys-Ing}, each of \eqref{SSPov} and \eqref{SSPext}, equivalent to \eqref{SSP}, gives a mereological analogue of the extensionality principle with respect to the relations $\Ov$ and $\Ext$, respectively:
\begin{gather*}
\Ko_{x,y\in U}\bigl(\Ko_{u\in U}(u\Ov x\Row u\Ov y)\IMP x=y\bigr), \label{ext-Ov} \tag{$\nz{ext}_{\scriptscriptstyle\pOv}$}\\
\Ko_{x,y\in U}\bigl(\Ko_{u\in U}(u\Ext x\Row u\Ext y)\IMP x=y\bigr). \label{ext-Ext} \tag{$\nz{ext}_{\scriptscriptstyle\pExt}$}
\end{gather*}
The above implications are invertible due to the property of `$=$'. It is proven that under the given definitions and assumptions $\PP$, both of the above principles are equivalent to \eqref{fun-Sum}. Furthermore, this shows once again that \eqref{SSP} entails both \eqref{ext-Ov} and \eqref{ext-Ext}, and that each of them entails \eqref{WSP}. Finally, notice that from \eqref{antys-Ing} and \eqref{PPP}, we get the principle \eqref{ext-PP}, which we already got from \eqref{fun-Sum} \citep[see][pp.~73 and 77, 86 and 94, 59, 47, respectively]{ja00b, ja18, ja13, ja20}.
\subsection{Featured theories\label{4.5}}
From the principles given in this section, 12 non-equivalent theories can be created, which have \eqref{nepe} as their thesis. Among them, only nine theory satisfy \eqref{EExt}. However, we believe that only those in which the principle \eqref{fun-Sum} holds deserve the name of ``theory of parts''. We have only three such theories. We will rank them from the weakest to the strongest. In each of them, the relation $\PP$ strictly partially orders universes, i.e., it is transitive and irreflexive and, therefore, also asymmetric. Therefore, we will define a given theory by listing additionally adopted assumptions:
\begin{enumerate}[1.]
\item \eqref{fun-Sum} \hfill equivalent to each of \eqref{ext-Ov} and \eqref{ext-Ext}
\item \eqref{fun-Sum}+\eqref{PPP}
\item \eqref{SSP} \hfill equivalent to each of ${\Sum}\subseteq{\Sup}$, \eqref{SSPov} and \eqref{SSPext}
\end{enumerate}
Even the weakest of them contains all the discussed formulas except \eqref{PPP}, \eqref{SSP} and their equivalents. Of course, among the three distinguished theories, the strongest one is the most interesting. The full grid (lattice) of theories built from the given formulas can be found in \citep[pp.~54, 42, respectively]{ja13, ja20}.

Let us note that all of the theories under consideration can be axiomatized in an elementary way, i.e., without using the predicate `$\in$' and variables for subsets of a universe of considerations. We replace \eqref{fun-Sum} with \eqref{ext-Ov} and remove notations of the form `$\ldots\in U$'.

However, a problem arises: Are these three theories really existentially neutral? We have an affirmative answer, even in the case of the strongest one, because the principle (SSP) is existentially neutral. Although it postulates the existence of some object, it is connected with the property of the relation of \emph{being a part of}. That the principle \eqref{SSP} does not postulate the existence of mereological sums is further evident from the fact that it is equivalent to \eqref{SSPov} and \eqref{SSPext}. These latter certainly do not postulate the existence of mereological sums.

Let us add that principle \eqref{SSP} also gives equivalence to two explications of the relational notion of \emph{being a collective class of}. Thus, the theory with \eqref{SSP} has the greatest merits in being recognized as the proper existentially neutral theory of parts.

\section{Existentially involved theories of parts\label{sec5}}
All of the existentially involved theories presented here will be stronger than the strongest of the existentially neutral theories, i.e., the principle \eqref{SSP} will hold in them. In the introduction, we mentioned that such theories may have assumptions that implicitly postulate the existence of collective classes (as mereological sums). Such will be the case for the first of the existentially involved theories presented here.

\subsection{Two theories of mereological strict partial orders}
Two theories of mereological strict partial orders. The first of them is obtained by adding the equality ${\Sum}={\Sup}$ to the axioms of the theory of strict partial orders, i.e., to \eqref{pz-PP} and \eqref{p-PP}. In Subsection~\ref{4.2}, we showed that the inclusion ${\Sum}\subseteq{\Sup}$ is equivalent to \eqref{SSP}. We also mentioned that it is the inclusion ${\Sup}\subseteq{\Sum}$ that makes us obtain an existentially committed theory. Indeed, the following example shows that the holding of this inclusion forces the existence of mereological sums that we do not obtain from the accepted definitions and principles \eqref{pz-PP}, \eqref{p-PP} and \eqref{SSP}.

Consider a universe that is a sharp partial order satisfying the principle \eqref{SSP} and has at least four elements, among which there is the unity $\mathbf{1}$ and three pairwise non-overlapping objects $o_1$, $o_2$ and $o_3$, which are furthermore only parts of the unity (i.e.\ lie directly under~$\mathbf{1}$). The unity $\mathbf{1}$ is the supremum of each of the pairs $\{o_1,o_2\}$, $\{o_1,o_3\}$, $\{o_2,o_3\}$, but none of them has a mereological sum. Hence, the inclusion ${\Sup}\subseteq{\Sum}$ does not hold. To save it, we need to introduce three new objects $\llbracket o_1,o_2\rrbracket$, $\llbracket o_1,o_3\rrbracket$, $\llbracket o_2,o_3\rrbracket$ which are the sums of pairs $\{o_1,o_2\}$, $\{o_1,o_3\}$ and $\{o_2,o_3\}$, respectively.

Let us also note that the inclusion ${\Sup}\subseteq{\Sum}$ entails the universe to be non-degenerate, i.e., it has at least two elements. Indeed, if $U=\{u\}$, then $u$ is the supremum of the empty set. Therefore, it would also have to be a mereological sum of this set, and we know that there is no such sum.

If we allow degenerate structures, then we have to restrict the inclusion of ${\Sup}\subseteq{\Sum}$ to non-empty sets, i.e. assume the following condition:
\[\label{eq:Sup=>Sum+}\tag{${\dag}$}
\Ko_{S\in\mathcal{P}(U)}\Ko_{x\in U}\bigl((S\ne\emptyset\wedge x\Sup S)\IMP x\Sum S\bigr).
\]
We thus obtain a theory equivalent to the theory of strict partial orders with the added condition which says that the relations $\Sum$ and $\Sup$ coincide on non-empty sets:
\[\label{eq:Sup=Sum+}\tag{${\ddag}$}
\Ko_{S\in\mathcal{P}(U)}\Ko_{x\in U}\bigl(x\Sum S\iff (S\ne\emptyset\wedge x\Sup S) \bigr).
\]
\indent In both cases, we obtain existentially involved theories that deserve to be called \emph{mereological theories of strict partial orders}.\footnote{Let us remember, however, that in both cases, the added conditions are not definitions of $\Sum$.}
\subsection{Simons' Minimal Extensional Mereology\label{sec:memS}}
In \citep[p.~31]{S87}, `Minimal Extensional Mereology' is the name given to the theory, which is based on axioms \eqref{pz-PP}, \eqref{p-PP}, \eqref{WSP} and the following:
\[\label{ax-B3}\tag{$\nz{c}{\exists}\mathord{\bprod}$}
\Ko_{x,y\in U}\bigl(x\Ov y\IMP\Ks_{z\in U}\Ko_{u\in U}(u\Ing z\Row u\Ing x \AND u\Ing y)\bigr),
\]
which states the conditional existence of an object, this being the ``product'' of two overlapping objects. Observe that \eqref{ax-B3} can be strengthened to the following equivalence:
\[
\Ko_{x,y\in U}\bigl(x\Ov y\iff\Ks_{z\in U}\Ko_{u\in U}(u\Ing z\Row u\Ing x\AND u\Ing y)\bigr).
\]
We have already mentioned that given \eqref{WSP}, it is unnecessary to take \eqref{pz-PP} as an axiom.

It has been proven that \eqref{SSP} is a thesis of Minimal Extensional Mereology \citep[see, e.g.,][pp.~120, 148, 77, 67, respectively]{ja00b, ja18, ja13, ja20}. Therefore, we have a thesis \eqref{fun-Sum} in it. Axiom \eqref{ax-B3} is related to the relation $\Sum$. It says that for any overlapping objects $x$ and $y$, the product $x\bprod y$ is the unique mereological sum of all their common ingredients, i.e., $x\bprod y$ is the sum of the non-empty set $\{u\in U : u\Ing x \AND u\Ing y\}$. Let us add that \eqref{ax-B3} cannot be obtained from \eqref{SSP} \citep[see, e.g.,][pp.~144, 183, 77, 68, respectively]{ja00b, ja18, ja13, ja20}.

It has also been shown that minimal extensional mereology crosses with both theories of mereological strict partial orders. More precisely, the former includes neither the inclusion ${\Sup}\subseteq{\Sum}$ nor condition \eqref{eq:Sup=>Sum+}. The latter does not include \eqref{ax-B3}. Therefore, one can consider stronger theories that arise by adding to the former either the inclusion ${\Sup}\subseteq{\Sum}$ or condition~\eqref{eq:Sup=>Sum+} \citep[see][pp.~147--149, 170--181, 78--80, 70, respectively]{ja00b, ja18, ja13, ja20}.

\subsection{Minimal Closure Mereology}
We have already mentioned in the introduction that existentially neutral theories do not postulate the existence of a collective set formed by two objects that are parts of a third. In the terminology we have adopted, such a postulate has the following form:
\[\label{cEs}\tag{$\nz{c}{\exists}{\bsum}$}	
\Ko_{x,y\in U}\bigl(\Ks_{u\in U}(x\Ing u\AND y\Ing u)\IMP\Ko_{z\in U}\;z\Sum \{x,y\}\bigr).
\]
Here, we have the conditional existence of a sum of objects being the ingredienses of a third object. We add this condition to Minimal Extensional Mereology, in which---as we remember---we have \eqref{fun-Sum}. Thus, the unique mereological sum of such objects $x$ and $y$ can be denoted by $x\bsum y$ or $\llbracket x,y\rrbracket$.

Minimal Closure Mereology was studied by \citet{S87} and \citet{CV}. There, however, a variant of the axiom \eqref{cEs} was adopted, in which the condition `$z\Sum\{x,y\}$' was replaced by an equivalent condition: $\Ko_{u\in U}(u\Ov z \Row(u\Ov x\OR u\Ov y))$. This equivalence is obtained from the general thesis \eqref{DolarOv} taken for $S=\{x,y\}$.

\subsection{Grzegorczykian mereology\label{subsec:mtG}}
Grzegorczyk presented his system of mereology in an article \citeyearpar{G55}. In short, we can say that his theory differs from Leśniewski's mereology in that the former postulates the existence of mereological sums only for a finite number of objects. In contrast, the latter postulates the existence of mereological sums for any number of objects, even if there is an infinite number of them. We will present Grzegorczyk's theory in a simplified form using the terminology we have adopted. However, this approach is equivalent to Grzegorczyk's original approach \citep[see, e.g.,][Subsection~III.7]{ja13, ja20}. Both approaches are elementary (in the sense of Subsection~\ref{4.5}).

We can assume that Grzegorczykian mereology is a theory of strict partial orders\footnote{In the original, the relation $\Ing$ was assumed as primary, i.e., these are partial orders.} in which instead of axiom \eqref{SSP}, its stronger version is adopted in the form of the Strong Super-Supplementation Principle:
\[\label{SSSP}\tag{$\nz{SSP+}$}
\Ko_{x,y\in U}\bigl(x\nIng y\IMP\Ks_{z\in U}\bigl(z\Ing x\AND  z\Ext y\AND\Ko_{u\in U}(u\Ing x\AND  u\Ext y\Imp u\Ing z)\bigr)\bigr).
\]
Obviously, the <<normal>> versions is logically entailed by the <<super>> version. Moreover, in \citep{ja00b, ja18, ja13, ja20}, it was proved that it also entails principle \eqref{eq:Sup=Sum+}, i.e., the relations $\Sum$ and $\Sup$ coincide on non-empty sets.

Moreover, in the theory under consideration, we have the uniqueness of mereological sum expressed by \eqref{fun-Sum}. Axiom \eqref{SSSP} is also related to the relation $\Sum$. It says that for arbitrary objects $x$ and $y$ such that $x$ is not an ingrediens of $y$, there is exactly one mereological sum of all ingredients of $x$ that are exterior to $y$, i.e., the mereological sum of the set $\{u \in U: u\Ing x \AND u\Ext y\}$ \citep[see][pp.~96, 90, respectively]{ja13, ja18}. This sum can therefore be considered as the mereological difference of $x$ and $y$, which we denote by $x - y$.

Grzegorczyk also adopted an axiom that is equivalent to the assumption of the unconditional existence of the mereological sum of any two objects:
\[\label{Es}\tag{${\exists}{\bsum}$}	
\Ko_{x,y\in U}\Ko_{z\in U}\;z\Sum \{x,y\}.
\]
Since we have \eqref{fun-Sum}, the unique mereological sum of arbitrary objects $x$ and $y$ can be denoted by $x\bsum y$ or $\llbracket x,y\rrbracket$. Hence we get that we also have mereological sums of any non-empty finite subset of the universe. Simply, using \eqref{Es}, we sum successive elements of a given non-empty finite subset of the universe.

Grzegorczyk adopted yet another axiom, which is equivalent to principle \eqref{ax-B3} from Minimal Extensional Mereology. However, as shown, the adoption of this axiom is unnecessary. Given the operation of mereological difference, we can obtain the mereological product of any overlapping objects $x$ and $y$ in the following way \citep[see][pp.~99 and~116, 93 and~110, respectively]{ja13, ja20}:
\[
x\bprod y = \begin{cases}
x &\text{if $x\Ing y$,}\\
y &\text{if $y\Ing x$,}\\
x- (x-y) &\text{if $y\POv x$.}
\end{cases}
\]
\indent From the above-adopted axioms we will not obtain the existence of unity. Therefore, we admit such structures in which there is no object encompassing all the remaining elements, as its parts. The class of models of Grzegorczykian mereology was compared with a certain class of lattices with zero, which in \citep{ja13, ja20} were called Grzegorczykian lattices. It was proved that this class of models of Grzegorczykian mereology coincides with the class of structures that arise from non-degenerate Grzegorczykian lattices after removing zero from them. It was also proved that conversely, from each model of Grzegorczykian mereology after adding a zero element to it we will obtain a non-degenerate Grzegorczyk lattice \citep[see][]{ja13, ja20}.

In the class of models of Grzegorczykian mereology, a subclass of structures with unity can be distinguished. It corresponds to a theory that arises by adding an axiom postulating the existence of unity:
\[
\Ks_{x\in U}\Ko_{u\in U}\;  u\Ing x.
\]
Let us denote the unity postulated in this way by $\mathbf{1}$. We know that it is the mereological Sum of the universe, i.e. we have $\mathbf{1}\Sum U$. Since in all models of Grzegorczykian mereology, we have a mereological difference; therefore, in models with unity, we also have a mereological complement of an object different from $\mathbf{1}$. Namely, for any such $x$, we put $-x := \mathbf{1}-x$.

\looseness=-1 The class of models of Grzegorczykian mereology with unity was compared with the class of Grzegorczykian lattices with unity. However, it was proved that the latter overlaps the class of Boolean lattices---equivalents of Boolean algebras \citep[see][]{ja13, ja20}.\footnote{Every Boolean lattice (Boolean algebra) has a unity and a zero.} It was proved that the class of models of Grzegorczykian mereology with unity overlaps the class of structures that arise from non-degenerate Boolean lattices after removing zero from them. It was also proved that conversely, from any model with unity, after adding zero to it, we obtain a non-degenerate Boolean lattice \citep[see][]{ja13, ja20}.

There is a fundamental difference between models without unity and models with unity. It has been shown that if a given structure does not have unity, then it cannot have it, i.e., it is impossible to add to it such an element that in the extended structure would be a unity. Indeed, for arbitrary $x,y\in U$, we have $x\bsum y\in U$. Therefore, if we add to the set $U$ some additional element $\mathit{1}$, which would be a unity, then for $x$, we will not find in $U$ such $y$ that $x\bsum y =\mathit{1}$, i.e., $x$ cannot have a complement in such an extended structure \citep[{cf.}][]{ja13,ja20}.

Of course, all finite Grzegorczykian mereological structures have a unity that is the sum of the entire universe.

\subsection{Classical mereological structures\label{subsec:KSM-ML}}
As we mentioned in the introduction, Leśniewskian mereology can be translated into the language of the theory of relational structures. In this form, we assume for it the axioms of strict partial orders, the uniqueness of mereological sum in the form \eqref{fun-Sum} and the axiom postulating the existence of the mereological sum for every non-empty subset of the universe:
\[\label{Esum}\tag{${\exists}{\Sum}$}
\Ko_{S\in\mathcal{P}(U)}(S\ne\emptyset\Imp\Ks_{x\in U}\;x\Sum S).
\]
We can assume the existence of a mereological sum only for non-empty sets.

The structures satisfying the axioms of strict partial orders,  \eqref{fun-Sum} and \eqref{Esum} are called \emph{classical mereological structures}, and their theory---\emph{classical mereology}. In every such structure, we have the unity $\mathbf{1}$ (where $\mathbf{1}\Sum U$) since the universe $U$ has the mereological sum. Moreover, all the principles considered so far are theses of classical mereology, so the relations $\Sum$ and $\Sup$ coincide on non-empty sets \citep[{cf.}][]{ja00b, ja13, ja18, ja20}. It is also known that classical mereology is not elementary axiomatizable \citep[see][]{ja00b}.

\citet{T29, T56} showed that classical mereology corresponds to the theory of complete non-degenerate Boolean lattices, i.e., in which every subset of the universe has a supremum. Namely, the class of all classical mereological structures coincides with the class of structures that arise from complete and non-degenerate Boolean lattices after removing zero from them. Conversely, from every classical mereological structure, after adding zero to it, we obtain a non-degenerate complete Boolean lattice \citep[{cf.}][]{ja00b, ja13, ja18, ja20}.

The only difference between Leśniewskian mereology and Grzegorczykian mereology with unity is that the former postulates the existence of a mereological sum also for non-empty infinite sets. It shows that Grzegorczykian finite mereological structures are the same as finite classical mereological structures.
\section{A problem related to the transitivity of the part-relation\label{sec6}}

In Subsection~\ref{sec:czjk}, we already mentioned the problems with the transitivity of the notion of \emph{being part of}. In the book \emph{Semantics}, John \citet{Ly77} likened the fact that an object $x$ is a part of an object $y$ to the semantic correctness of a sentence of the form ‘$y$ has $x$’. For example, the following sentences are semantically correct:
\begin{enumerate}[\textbullet]
\item The orchestra ($z$) has a violin section ($x$).\hfill  $x\PP z$
\item The orchestra ($z$) has a violinist ($y$).\hfill  $y\PP z$
\item The violinist has a heart ($u$).\hfill  $u\PP y$
\item The violinist has an arm ($v$).\hfill  $v\PP y$
\end{enumerate}
The following forms are not, however, correct:
\begin{enumerate}[\textbullet]
\item The orchestra has the arm of the violinist.
\item The violin section has the violinist's heart.
\end{enumerate}

In cases when is problematic the transitivity of the relation \emph{being part of}, we propose to assume that $\PP$ is acyclic, i.e., satisfies \emph{ac-PP}, and that it is locally transitive in the following sense. If $x$ is part of $y$, then the transitivity is to hold on any path leading from $x$ to $y$, which consists of objects that are parts of successive ones appearing on this path. In other words, if we have $x\PP y$ and $x\PP z_1 \PP z_2 \PP \cdots \PP z_n\PP y$, then the relation $\PP$ is transitive on the set $\{x,z_1,\ldots, z_n,y\}$. By virtue of the acyclicity of the relation $\PP$, no such path is closed, i.e. does not lead from a given object to itself.

For example, let $y$ be a violinist, $x$ be a finger of his right hand,  $z_1$ be his right hand, and $z_2$ be his right arm. Then $x\PP y$, $x\PP z_1\PP z_2\PP y$, $x\PP z_2$ and $z_1\PP y$. On the other hand, let $y$ be an orchestra, $x$ be a violinist in the first violin section of this orchestra, $z_1$ be the first violin section in this orchestra, and $z_2$ be the string section in this orchestra. Then $x\PP y$, $x\PP z_1\PP z_2\PP y$, $x\PP z_2$ and $z_1\PP y$. An orchestra is a system of parts that make a direct functional contribution to the whole. In this system, the musicians and the conductor are disjoint elements that are minimal with respect to the relation of \emph{being a part of}. Each member of a given orchestra is also such a system. Similarly, to use Rescher's example, a cell is also such a system of parts. Among these parts is its nucleus.

If we do not assume the transitivity of the relation of \emph{being a part of}, then in addition to its acyclicity and local transitivity, we also assume other axioms, among them those concerning maximally closed sets with respect to this relation. If we assume the transitivity of this relation, the entire universe is the only maximally closed set with respect to this relation.

In Chapter~4 of the books \citep{ja13, ja20} and in the article \citep{ja14}, various possible solutions to the construction of the theory of parthood without assumed transitivity and with assumed local transitivity were given.

\renewcommand{\bibfont}{\small}
\setlength{\bibsep}{5pt}

\end{document}